\documentclass[12pt,fleqn]{article}

\usepackage{geometry}

\usepackage[utf8]{inputenc}
\usepackage[T2A]{fontenc}
\usepackage[english,russian]{babel}

\usepackage{amsmath}
\usepackage{amssymb}
\usepackage{amsthm}
\usepackage{array}
\usepackage[auth-sc]{authblk}
\usepackage{bm}
\usepackage{cite}
\usepackage{color}
\usepackage{dsfont}
\usepackage{enumerate}
\usepackage{enumitem}
\usepackage{epsf}
\usepackage{euscript}
\usepackage{graphicx} 
\usepackage{indentfirst}
\usepackage{latexsym}
\usepackage{mathrsfs}
\usepackage{mathtools}
\usepackage{tikz}
\usepackage{pgf}

\newtheorem{Prop}       {Предложение}

\newtheorem{Cor}        {Следствие}  
\newtheorem{Rem}        {Замечание}

\begin{document}

\selectlanguage{russian}

\title{Трюк Крипке и разрешимость монадических фрагментов модальных и суперинтуиционистских предикатных логик}

\author[1]{М.\,Рыбаков}
\author[2]{Д.\,Шкатов}
\affil[1]{ИППИ РАН, НИУ ВШЭ, ТвГУ}
\affil[2]{University of the Witwatersrand, Johannesburg} 

\date{}

\maketitle


\begin{abstract}
Трюк Крипке позволяет моделировать бинарную предикатную букву в классических формулах модальными формулами с двумя унарными предикатными буквами. Рассматриваются вариации трюка Крипке и возможности его применения в модальных и суперинтуиционистских предикатных логиках. Кроме того, обсуждаются ситуации, когда применить трюк Крипке невозможно.
\end{abstract}

Сол Крипке заметил~\cite{Kripke62}, что в классических формулах первого порядка бинарная буква моделируется с помощью двух унарных: достаточно заменить формулы вида $P(x,y)$ формулами вида $\Diamond(Q_1(x)\wedge Q_2(y))$. Такое моделирование позволяет погрузить классическую логику бинарного предиката в монадический фрагмент любой модальной логики, расширяющей классическую логику предикатов и содержащуюся в~$\mathbf{QS5}$. Поскольку классическая логика бинарного предиката неразрешима, получаем неразрешимость монадических фрагментов большого класса модальных логик. Учитывая неразрешимость теории симметричного иррфелексивного предиката~\cite{Sper:2016,MR:2022DAN}, использование формулы $\neg\Diamond(Q(x)\wedge Q(y))$ вместо $\Diamond(Q_1(x)\wedge Q_2(y))$ даёт неразрешимость модальных логик одного унарного предиката, содержащихся в $\mathbf{QS5}$, причём лишь при трёх предметных переменных в языке.

Анализ трюка Крипке показывает, что возможность его применения предполагает выполнение некоторых условий, в частности:
\begin{itemize}[leftmargin=5em]
\item[$(\mathit{TK\!}_1)$]
использование формул, где под модальностью может быть более одной свободной переменной;
\item[$(\mathit{TK\!}_2)$]
отсутствие в логике формул, ограничивающих число миров, достижимых из произвольного мира.
\end{itemize}

Нарушение условия $(\mathit{TK\!}_1)$ приводит к т.н.\ мон\b{о}дическим фрагментам (когда в области действия модальности могут находится формулы с не более чем одной свободной переменной), которые часто оказываются разрешимыми~\cite{WZ01}. 

Как заметили авторы, нарушение условия $(\mathit{TK\!}_2)$ приводит к разрешимости мон\b{а}дических фрагментов (когда в формулах допускаются предикатные буквы валентности не более чем один)~\cite{MR:2017LI} и даже с равенством~\cite{RSh:2023JLC}. Отметим, что в семантике Крипке для модальных и суперинтуиционистских логик равенство при этом можно понимать по-разному, и, в отличие от классической логики предикатов, эти понимания приводят к разным множествам истинных формул. Так, в шкалах Крипке (с~предметными областями) равенство может определяться~\cite{GShS} как
\begin{itemize}[leftmargin=5em]
\item[$(\mathit{Eq}_1)$]
наследственная вверх конгруэнтноть, 
\item[$(\mathit{Eq}_2)$]
наследственная вверх и вниз конгруэнтность, 
\item[$(\mathit{Eq}_3)$]
предикат совпадения. 
\end{itemize}
Принципы $(\mathit{Eq}_1)$, $(\mathit{Eq}_2)$ и $(\mathit{Eq}_3)$ дают семантики, различимые как модальными, так и интуиционистскими формулами: существует шкала Крипке, для которой эти принципы приводят к трём разным множествам истинных в ней интуиционистских (а следовательно, и модальных) формул~\cite{GShS}. Принцип $(\mathit{Eq}_1)$ выражается модальной формулой $(x=y)\to \Box(x=y)$, а принцип $(\mathit{Eq}_2)$~--- формулой $(x=y)\leftrightarrow \Box(x=y)$. В интуиционистской семантике принцип $(\mathit{Eq}_1)$ выполняется автоматически, а принцип $(\mathit{Eq}_2)$ описывается формулой $(x=y)\vee \neg(x=y)$ и соответствует понятию разрешимого равенства. Принцип~$(\mathit{Eq}_3)$ не является определимым ни в одном из этих языков. Тем не менее, для многих естественных классов шкал Крипке принципы $(\mathit{Eq}_2)$ и $(\mathit{Eq}_3)$ не различимы ни модально, ни интуиционистски~\cite{GShS}. Приведённые ниже утверждения справедливы для каждого из этих трёх пониманий равенства.

\begin{Prop}
\label{prop:1}
Монадический фрагмент с равенством модальной логики шкалы Крипке с конечным множеством миров алгоритмически разрешим.
\end{Prop}

\begin{Cor}
Монадические фрагменты с равенством модальных и суперинтуиционистских логик рекурсивно перечислимых классов шкал с конечным числом миров принадлежат классу\/ $\Pi^0_1$, т.е. имеют рекурсивно перечислимое дополнение.
\end{Cor}

Учитывая~\cite{RShJLC20a,RShJLC21b}, получаем ещё одно следствие.
Если $L$~--- модальная или суперинтуиционистская логика, то $L_{\mathit{wfin}}^=$ обозначает, соответственно, модальную или суперинтуиционистскую логику с равенством, определяемую классом шкал логики $L$, содержащих конечное множество миров; при этом семантически равенство определяется в соответсвии с любым из принципов $(\mathit{Eq}_1)$, $(\mathit{Eq}_2)$ или $(\mathit{Eq}_3)$.

\begin{Cor}
Пусть $L$~--- одна из модальных логик\/
$\mathbf{QK}$, $\mathbf{QT}$, $\mathbf{QK4}$, $\mathbf{QK5}$, $\mathbf{QK45}$, $\mathbf{QKD}$, $\mathbf{QKD4}$, $\mathbf{QKD45}$, $\mathbf{QKB}$, $\mathbf{QKTB}$, $\mathbf{QS4}$, $\mathbf{QS5}$, $\mathbf{QGL}$, $\mathbf{QGrz}$ или одна из супеинтуиционистских логик\/ $\mathbf{QInt}$, $\mathbf{QKP}$, $\mathbf{QKC}$. Тогда монадический фрагмент логики $L_{\mathit{wfin}}^=$~--- как с равенством, так и без равенства~--- является\/ $\Pi^0_1$-полным.
\end{Cor}

Учитывая полноту по Крипке логик вида $\mathbf{QAlt}_n^=$~\cite{ShSh:2023JLC}, получаем следующее следствие.

\begin{Cor}
Для любого $n\in\mathds{N}$ монадический фрагмент с равенством модальной логики $\mathbf{QAlt}_n^=$ является алгоритмически разрешимым.
\end{Cor}

\begin{Rem}
Приведённые следствия останутся справедливым, если рассматривать семантику постоянных областей.
\end{Rem}

Тем не менее, некоторые ослабления приведённых условий всё-таки позволяют применить трюк Крипке (или какую-то его модификацию). Например, в шкалах Крипке с конечным числом миров действует второе из указанных ограничений; но если при этом для каждого $n\in\mathds{N}$ допускается наличие шкалы с миром, видящим не менее $n$ миров, то логика унарного предиката для такого класса шкал не будет не только разрешимой, но и рекурсивно перечислимой~\cite{Skv95,RShJLC20a,RShJLC21b}.

В интуиционистской логике применение трюка Крипке затруднительно, что связано с наличием отрицания. Тем не менее, в очень многих случаях отрицание можно промоделировать импликацией к новой пропозициональной букве, а в получившихся позитивных формулах использовать следующую модификацию трюка Крипке: формулы вида $P(x,y)$ заменить формулами вида $(Q_1(x)\wedge Q_2(y)\to p)\vee q$~\cite{RSh19SL,RShJLC21b,RShsubmitted}. Отметим, что такое моделирование позволяет остаться в рамках позитивного фрагмента; тем не менее, здесь пропозициональную букву $p$ можно заменить формулой~$\bot$, т.е.\ моделировать формулу $P(x,y)$ формулой $\neg(Q_1(x)\wedge Q_2(y))\vee q$.

Приведём некоторые результаты, полученные авторами для модальных и суперинтуиционистских логик с использованием различных модификации трюка Крипке:
\begin{itemize}
\item
многие модальные и суперинтуиционистские логики унарного предиката неразрешимы в языке с двумя предметными переменными~\cite{RSh19SL};
\item
многие модальные и суперинтуиционистские логики одного унарного предиката естественных классов конечных шкал Крипке не являются рекурсивно перечислимыми при трёх предметных переменных~\cite{RShJLC20a,RShJLC21b};
\item
модальные логики различных классов линейных шкал (сильно) неразрешимы в языке с одной унарной буквой, одной пропозициональной буквой и двумя предметными переменными~\cite{RSh20AiML,RShJLC21c,RSh21SR};
\item
модальные логики нётеровых порядков сильно неразрешимы в языке с двумя унарными буквами, одной пропозициональной буквой и двумя предметными переменными~\cite{Rybakov02,RybIGPL22}.
\item
полимодальные темпоральные логики унарного предиката сильно неразрешимы~\cite{MREK:2017SR} в языке с двумя предметными переменными~\cite{RShJLC21d}.
\end{itemize}
Кроме того, авторами готовится работа~\cite{RShsubmitted}, где предполагается описать особенности, связанные с использованием трюка Крипке и его модификаций в различных неклассических логиках.

\medskip

\textit{
Работа выполнена в ИППИ РАН при поддержке Российского научного фонда, грант \mbox{21--18--00195}.
}


\end{document}